\documentclass{amsart}

\usepackage{graphicx}

\usepackage[utf8]{inputenc}
\usepackage{amsmath,amsfonts,amssymb}
\usepackage{hyperref}
\usepackage{bm}
\usepackage{booktabs}
\usepackage{mathtools}
\usepackage{mathdots}
\usepackage{xcolor}
\usepackage{algorithm}
\usepackage{algpseudocode}
\usepackage{url}
\usepackage{pgfplots}
\pgfplotsset{compat=newest}
\usetikzlibrary{plotmarks}
\usetikzlibrary{arrows.meta}
\usetikzlibrary{decorations.markings}
\usepgfplotslibrary{patchplots}
\usepackage{grffile}
\definecolor{mycolor1}{rgb}{0.00000,0.44700,0.74100}%
\definecolor{mycolor2}{rgb}{0.85000,0.32500,0.09800}%
\definecolor{mycolor3}{rgb}{0.92900,0.69400,0.12500}%

\newcommand{\ri}{\mathrm{i}}
\newcommand{\e}{\mathrm{e}}
\newcommand{\R}{\mathbb{R}}
\newcommand{\C}{\mathbb{C}}

\newcommand{\cO}{\mathcal{O}}

\renewcommand{\tilde}{\widetilde}

\DeclareMathOperator{\diag}{diag}

\newtheorem{theorem}{Theorem}
\newtheorem{lemma}[theorem]{Lemma}

\newtheorem{remark}[theorem]{Remark}

\title[Fast optimization for data-driven predictive control]{Fast and memory-efficient optimization for large-scale data-driven predictive control}

\newcommand*{\mailto}[1]{\href{mailto:#1}{\nolinkurl{#1}}}

\author[P. Schmitz]{Philipp Schmitz}
\address{Technische Universität Ilmenau}
\email{\mailto{philipp.schmitz@tu-ilmenau.de}}
\urladdr{\url{https://www.tu-ilmenau.de/obc/team/philipp-schmitz}}

\author[M. Schaller]{Manuel Schaller}
\address{Technische Universität Ilmenau}
\email{\mailto{manuel.schaller@tu-ilmenau.de}}
\urladdr{\url{https://www.tu-ilmenau.de/deq/manuel-schaller}}

\author[M. Voigt]{Matthias Voigt}
\address{UniDistance Suisse}
\email{\mailto{matthias.voigt@fernuni.ch}}
\urladdr{\url{https://fernuni.ch/profil/matthias-voigt}}

\author[K. Worthmann]{Karl Worthmann}
\address{Technische Universität Ilmenau}
\email{\mailto{karl.worthmann@tu-ilmenau.de}}
\urladdr{\url{https://www.tu-ilmenau.de/obc/team/karl-worthmann}}


\begin{document}

\begin{abstract}
Recently, data-enabled predictive control (DeePC) schemes based on Willems' fundamental lemma have attracted considerable attention. At the core are computations using Hankel-like matrices and their connection to the concept of persistency of excitation.
We propose an iterative solver for the underlying data-driven optimal control problems resulting from linear discrete-time systems. To this end, we apply factorizations based on the discrete Fourier transform of the Hankel-like matrices, which enable fast and memory-efficient computations. 
To take advantage of this factorization in an optimal control solver and to reduce the effect of inherent bad conditioning of the Hankel-like matrices, we propose an augmented Lagrangian lBFGS-method. We illustrate the performance of our method by means of a numerical study.
\end{abstract}

\maketitle

\section{Introduction}

Data-driven methods have recently received a lot of attention in the field of optimal and predictive control, see, e.g., \cite{HewiWabe20,bold2023practical} as well as the survey articles~\cite{MartScho23,verheijen2023handbook}. Far-reaching results for linear discrete-time systems have been achieved using a key result, also coined as Willems' fundamental lemma and deduced by J.C.\ Willems, P.\ Rapisarda, I.\ Markovsky, and  B.L.M.\ De Moor in~\cite{willems2005note}. This is due to 
an exact representation of the input-output behavior by means of Hankel-like matrices, see, e.g., \cite{CoulLyge19,berberich2020data,schmitz2022willems,van2020willems} and the recent survey articles~\cite{markovsky2021behavioral,FaulOu23} and the references therein.
Herein, the input-output behavior is parameterized by suitably constructed Hankel-like matrices. Then, large-scale optimal control problems with highly structured matrices have to be solved to determine the sought-after optimal controls. However, the dimension of the problem and condition number rapidly grow with increasing number of data points such that memory- and time-efficient algorithms are required for a fast solution, e.g., in model predictive control. 
To this end, various approaches to reduce computational complexity and storage memory demand have been proposed for Hankel-based data-driven optimal control. One possibility is a compression of the data-matrix via singular value decomposition \cite{zhang2023dimension,yin2021maximumb}. Another approach proposes a segmentation into multiple shorter prediction horizons \cite{o2022data}, which leads to smaller Hankel-like matrices.

In this work, we show how one can alleviate the high dimensionality and thus the immense memory requirements by means of discrete Fourier transform (DFT)-based matrix-vector multiplications combined with an iterative optimal control solver. We leverage a factorization based on circulant matrices exploiting the Hankel structure of the underlying matrices in data-driven optimal control. 
Consequently, DFT allows to implement evaluations of matrix-vector products without setting up the system matrix explicitly and hence serves as a memory-efficient and fast foundation for iterative solvers. 

Since the convergence rate of iterative solvers often hinges on a good condition number of the system, preconditioned methods are necessary. As in addition to the bad conditioning of Hankel-based problems, the dimension is very large, preconditioning techniques that circumvent storing the matrix explicitly are necessary.
To tackle the high dimensionality and the bad conditioning, we propose an augmented Lagrangian-based approach with a limited-memory Broyden--Fletcher--Golfarb--Shanno (lBFGS) algorithm as inner solver. Both components of this solver are highly successful in large-scale optimal control, see, e.g., \cite{NoceWrig06,LiuNoce89}. The lBFGS-algorithm allows for an iterative build-up of a preconditioner without having to store the matrices explicitly and thus benefits from exploiting our efficient DFT-based matrix-vector multiplication.

The contribution of this work can thus be summarized as follows. First, we formulate a DFT-based factorization for fast and memory-efficient implementation of matrix-vector multiplication involving Hankel-like matrices. We then apply this factorization to data-driven optimal control. 
After formulating the proposed Lagrangian approach with lBFGS as inner solver, we provide a numerical study showcasing the convergence properties. Whereas other iterative solvers, such as the gradient method or MINRES, fail due to the bad conditioning, we illustrate that the proposed augmented Lagrangian lBFGS-method yields fast convergence while being memory-efficient.

\noindent \textbf{Notation.} For a finite sequence in $\R^d$ of length $N$ we use the short-hand notation $\bm u_N = (u_0,\dots, u_{N-1})$. The identity matrix in $\R^{d\times d}$ and the zero matrix in $\R^{m\times n}$ are denoted by $I_d$ and $0_{m\times n}$, respectively (where we sometimes omit subscripts if clear from context). For the Kronecker product of two matrices $A$ and $B$ we write $A\otimes B$.

\section{Motivation: Data-driven predictive control}
Recently, Willems et al.'s fundamental lemma \cite{willems2005note} has attracted a lot of interest in the data-driven control community as it paves the way to data-enabled controller design for discrete-time linear time-invariant systems \cite{CoulLyge19} of the form
\begin{equation}
\label{eq:LTI}
    x_{k+1} = Ax_k + Bu_k, \quad k=0,1,\ldots,
\end{equation}
with matrices $A\in\mathbb R^{n\times n}$ and $B\in\mathbb R^{n\times m}$ and where $x_k$ and $u_k$ denote the states and inputs, respectively. Under the assumption of controllability, the fundamental lemma allows for a complete, non-parametric, and data-based description of the system's behavior only by measured data including a persistently exciting input signal. An input sequence $\bm u_N$ is called \emph{persistently exciting of order $L$}, $L\geq 1$, if the Hankel-like matrix
\begin{equation*}
    H_L(\bm u_N) := \begin{bmatrix}
        u_0 & u_1& \dots & u_{N-L}\\
        u_1 & u_2& \dots & u_{N-L+1}\\
        \vdots & \vdots &  &\vdots\\
        u_{L-1} & u_L& \dots & u_{N-1}
\end{bmatrix}\in\mathbb R^{mL \times (N-L+1)}
\end{equation*}
has full row rank $mL$. Next we recall Willems et al.'s fundamental lemma, see~\cite{willems2005note}.
\begin{lemma}
\label{lem:fl}
    Suppose that the system~\eqref{eq:LTI} is controllable and let $(\tilde{\bm x}_N, \tilde{\bm u}_N)$ be a trajectory of~\eqref{eq:LTI} such that $\tilde{\bm u}_N$ is persistently exciting of order $L+n$. Then $(\bm x_L,\bm u_L)$ is a length-$L$ trajectory of system~\eqref{eq:LTI} if and only if there is $z\in\mathbb R^{N-L+1}$ satisfying
    \begin{equation}
    \label{eq:Hankel_description}
        \begin{bmatrix}
            x_0\\\vdots \\ x_{L-1}
        \end{bmatrix} = H_L(\tilde{\bm x}_N) z,\quad
        \begin{bmatrix}
            u_0\\\vdots \\ u_{L-1}
        \end{bmatrix} = H_L(\tilde{\bm u}_N) z.
    \end{equation}
\end{lemma}

A persistency of excitation order~$L+n$ results in a signal length 
\begin{align}\label{rm:length}
     N \geq (m+1)(L+n)-1.
\end{align}

A commonly used approach to control systems of the form~\eqref{eq:LTI} is model predictive control (MPC), where the feedback is computed by solving the optimal control problem (OCP)
\begin{equation}
\label{eq:ocp}
    \operatorname*{minimize}_{\bm x_L,\bm u_L}\ \frac{1}{2}\sum_{k = 0}^{L-1} \Big(x_k^\top Qx_k + u_k^\top Ru_k\Big)
\end{equation}
subject to the dynamics \eqref{eq:LTI} and some initial constraint $x_0 = x^0$ is solved. The matrices $Q$ and $R$ in the quadratic stage cost are assumed to be symmetric, positive definite. The fundamental lemma leverages the MPC approach to a data-enabled predictive control scheme, where the model described by the dynamics~\eqref{eq:LTI} is replaced by the description of 
 the finite-length trajectories in~\eqref{eq:Hankel_description} involving only data stored in Hankel-like matrices. In particular, there is no need to identify the system matrices $A$ and $B$. Given a data trajectory $(\tilde{\bm x}_N, \tilde{\bm u}_N)$ with an input signal $\tilde{\bm u}_N$ of sufficient persistency of excitation order, the OCP~\eqref{eq:ocp} can be equivalently formulated as the quadratic optimization problem
\begin{equation}
\label{eq:docp}
    \operatorname*{minimize}_z\ \frac{1}{2}z^\top \mathcal Sz\qquad \text{subject to } \mathcal P z = x^0,
\end{equation}
where
\begin{equation}
\label{eq:S}
    \mathcal S = \begin{bmatrix}
        H_L(\tilde{\bm x}_N)\\ H_L(\tilde{\bm u}_N)
    \end{bmatrix}^\top {\begin{bmatrix}
        I_{L}\otimes Q\\
        & I_{L}\otimes R
    \end{bmatrix}}\begin{bmatrix}
        H_L(\tilde{\bm x}_N)\\ H_L(\tilde{\bm u}_N)
    \end{bmatrix}
\end{equation}
and
\begin{equation*}
    \mathcal P = \begin{bmatrix}
    I_n & 0_{n\times n(L-1)}
\end{bmatrix}H_L(\tilde{\bm x}_N).
\end{equation*}

We emphasize that OCP~\eqref{eq:ocp} can be simply modified regarding setpoint tracking, which results in an additional linear term $q^\top z$ in the quadratic problem~\eqref{eq:docp}. As a consequence of the persistency of excitation, the matrix $\mathcal S$ has rank $mL + n$, while the rank of $\mathcal P$ equals $n$, cf.~\cite{willems2005note}.

Accounting for the required persistency of excitation order of at least $L+n$, the square matrix $\mathcal S$ has dimension $N-L+1\geq (L+n) m+n$, cf.~\eqref{rm:length}. Therefore, the memory demand for storing $\mathcal S$ as a dense numerical matrix, see Table~\ref{tab:mem_cons}, is proportional to \begin{equation*}
    (N-L+1)^2\geq\big((L+n) m+n\big)^2.
\end{equation*}

Depending on the available hardware this may render problem \eqref{eq:docp} infeasible to standard solvers (e.g.\ \texttt{MATLAB quadprog}) involving 'naive' matrix-vector multiplication for large system order $n$, high input dimension $m$, or long prediction horizon $L$. In this situation, iterative methods which allow for a matrix-free implementation are favorable.

\begin{table}[htb]
\label{tab:mem_cons}
    \centering
        \caption{The memory consumption of $\mathcal S$ stored as a dense numerical matrix and the trajectory $(\tilde{\bm x}_N ,\tilde{\bm u}_N)$ taking 8 bytes per entry for various problem sizes.}
    \label{tab:my_label}
    \begin{tabular}{rrrrr}
    \toprule
    $n$ & $m$ & $L$ & mem.\ $\mathcal S$ & mem.\ $(\tilde{\bm x}_N ,\tilde{\bm u}_N)$\\
    \midrule
        100 & 50 & 50 & 0.46 GB & 9.12 MB\\
        200 & 50 & 50 & 1.29 GB & 25.40 MB\\
        300 & 50 & 50 & 2.53 GB & 49.84 MB\\
        500 & 50 & 50 & 6.27 GB & 123.20 MB\\
        \midrule
        100 & 100 & 100 & 3.23 GB & 32.16 MB\\
        200 & 100 & 100 & 7.30 GB & 72.48 MB\\
        300 & 100 & 100 & 12.99 GB & 128.96 MB\\
        500 & 100 & 100 & 29.28 GB & 290.40 MB\\
        \bottomrule\\
    \end{tabular}
\end{table}

\section{DFT-based factorization of Hankel-like matrices}\label{sec:dft}
Let $\bm w_N$ denote the combined signal of $\tilde{\bm x}_N$ and $\tilde{\bm u}_N$ in $\mathbb R^d$, i.e., $w_k := \begin{bmatrix}
    \tilde x_k^\top & \tilde u_k^\top
\end{bmatrix}{}^\top$ and $d:=n+m$. A straight-forward approach to calculate the matrix-vector product $H_L(\bm w_N) v$ without storing the full Hankel-like matrix is to make use of the structure of $H_L(\bm w_N)$ and access its coefficients directly from the generating sequence $\bm w_N$. In the following we recall a factorization result for $H_L(\bm w_N)$, which solves the memory issue likewise and in addition is favorable for the speed of the matrix-vector multiplication.

\subsection{Factorization}
A row permutation of $H_L(\bm w_N)$ leads to the Toeplitz-like matrix
\begin{equation}
\label{eq:perm}
     T := \begin{bmatrix}
         w_{L-1} &  \dots & w_{N-1}\\
         \vdots &  & \vdots\\
         w_{0} & \dots & w_{N-L}
     \end{bmatrix}=\begin{bmatrix}
         &&I_d\\
         &\scalebox{-1}[1]{$\ddots$}\\
         I_d
     \end{bmatrix}H_L(\bm w_N),
\end{equation} which can be embedded into the circulant block-Toeplitz matrix $C\in \R^{Nd \times N}$,
\begin{equation*}
    C=\left[\begin{array}{cccc|cccc}
        w_{L-1} & w_L & \ldots & w_{N-1} & w_0 & w_1 & \ldots & w_{L-2} \\
         w_{L-2} & w_{L-1} & \ldots & w_{N-2} & w_{N-1} & w_0 & \ldots & w_{L-3} \\
        \vdots & \vdots & & \vdots & \vdots & \vdots & & \vdots \\
        w_{0} & w_1 & \ldots & w_{N-L} & w_{N-L+1} & w_{N-L+2} & \ldots & w_{N-1} \\ \hline
        w_{N-1} & w_0 & \ldots & w_{N-L-1} & w_{N-L} & w_{N-L+1} & \ldots & w_{N-2} \\
        w_{N-2} & w_{N-1} & \ldots & w_{N-L-2} & w_{N-L-1} & w_{N-L} & \ldots & w_{N-3} \\ 
        \vdots & \vdots & & \vdots & \vdots & \vdots & & \vdots \\
        w_L & w_{L+1} & \ldots & w_0 & w_1 & w_2 & \ldots & w_{L-1}
       \end{array}\right],
\end{equation*}
that is, $T=\begin{bmatrix}
    I_{Ld} & 0
\end{bmatrix}C \begin{bmatrix}
    I_{N-L+1} & 0
\end{bmatrix}^\top$. The matrix $C$ can be factorized as 
\begin{align*}
    C = \tilde{F} \Lambda F^{-1}
\end{align*}
where $F:=\big[\e^{-2\pi \ri kj/N}\big]_{k,j = 0,\dots, N-1}$ is a Fourier matrix in $\R^{N\times N}$, $\tilde F\ := F\otimes I_d\in\C^{Nd \times Nd}$ and $\Lambda := \diag(\Lambda_0,\,\ldots,\,\Lambda_{N-1}) \in \C^{Nd \times N}$ with
\begin{equation*}
 \begin{bmatrix} \Lambda_0 \\ \vdots \\ \Lambda_{N-1}  \end{bmatrix} := \tilde{F} \begin{bmatrix} w_{L-1}\\ \vdots \\ w_{N-1} \\ w_0 \\ \vdots \\ w_{L-2}  \end{bmatrix}, 
\end{equation*}
see also \cite[Section 4.8]{golub2013}.
\subsection{Complexity}
The matrix-vector product $Cv$ can be efficiently computed, since $F^{-1}v = \frac{1}{N}F^*v$ can be evaluated in $\cO(N \log N)$
time, $\Lambda v$ can be evaluated in $\cO(Nd)$ time and 
\begin{align*}
\tilde{F}v = (F \otimes I_d)v =  \left[\left(\e^{-2\pi\ri/N} \right)^{kj} I_d\right]_{k,j=0,\ldots,N-1} v 
\end{align*}
can be evaluated in $\cO(Nd \log N)$ time. As a consequence, the matrix-vector product $H_L(\bm w_N)z$ can be efficiently calculated via the following three steps: 1.)\ zero-padding $v=\begin{bmatrix}
    z^\top& 0
\end{bmatrix}{}^\top$, 2.)\ computing $Cv$, 3.)\ truncating the result and permuting its components, cf.~\eqref{eq:perm}. Therefore, the matrix-vector multiplication $H_L(\bm w_N)w$ has complexity $\cO(Nd\log N)$. The same complexity applies to matrix-vector products involving $H_L(\bm w_N)^\top$ and $\mathcal S$ in \eqref{eq:docp}.

In terms of the parameters $n$, $m$, $L$, this complexity can be also be expressed as
\begin{equation*}
    \cO\Big(\big((L+n)m+n\big)(n+m)\log\big((L+n)m+n\big)\Big).
\end{equation*}

In contrast to that, naive matrix-vector multiplication $\mathcal S z$ using the factorization in~\eqref{eq:S} has complexity 
\begin{equation*}
    \cO(L(n+m)(N-L+1)) = \cO(L(n+m)(m(L+n)+n)),
\end{equation*} while a precomputed matrix $\mathcal S$ leads to the complexity 
\begin{equation*}
    \cO((N-L+1)^2)=\cO((m(L+n)+n)^2).
\end{equation*}
In Table~\ref{tab:2} we list the complexities in the cases where one of the individual parameters $n$, $m$, or $L$ tends to infinity.
\begin{table}[htb]
    \centering
        \caption{Comparison of the complexity of DFT-based and naive matrix-vector multiplication $\mathcal Sz$ w.r.t.\ the state dimension~$n$, the input dimension~$m$, and prediction horizon $L$, where in each case the remaining two parameters are fixed.}
    \label{tab:2}
    \begin{tabular}{llll}
    \toprule
    & \multicolumn{3}{c}{complexity}\\
    \cmidrule{2-4}
    method & $n\rightarrow \infty$ & $m\rightarrow \infty$ & $L\rightarrow \infty$\\
    \midrule
    DFT-based & $\cO(n^2\log n)$ & $\cO(m^2\log m)$ & $\cO(L\log L)$\\
    naive, using~\eqref{eq:S} & $\cO(n^{2})$ & $\cO(m^2)$ & $\cO(L^2)$\\
    naive, precomp.\ $\mathcal S$ & $\cO(n^2)$ & $\cO(m^2)$ & $\cO(L^2)$\\
    \bottomrule\\
    \end{tabular}
\end{table}

\begin{remark}
\label{rm:fft}
    Most standard implementations of the DFT run significantly faster when the signal length $N$ has only small prime factors (not greater than $7$). Therefore, apart from the persistency of excitation this fact should be taken into account while choosing the length of the data trajectory.
\end{remark}
\section{Augmented Lagrangian with lBFGS}
Here, we briefly describe the solution method which enables a memory-efficient algorithm for solving the equality constrained optimization problem \eqref{eq:docp}. As an outer loop, we implement an augmented Lagrangian method, where in an inner loop, the arising unconstrained problems of the outer loop are solved by the Broyden--Fletcher--Goldfarb--Shanno (BFGS) method. We highlight that the combination of these two algorithms allows for solving \eqref{eq:docp} without storing the matrix and by just implementing a function that evaluates the matrix-vector products based on the given data and DFT. This is due to the fact that, whereas the BFGS method is a quasi-Newton method and hence does not suffer from ill-conditioning as much as the gradient method, the preconditioner is built up successively via rank-one updates. These updates of the preconditioner, given by an approximation of the inverse Hessian, can be performed efficiently with the DFT-based factorizations presented before.

The presentation in this section is based on \cite{NoceWrig06}. In this part, we abbreviate $q := N-L+1$ being the dimension of the variable $z$ in \eqref{eq:docp}.

\subsection{Outer loop: Augmented Lagrangian}
To solve the linearly constrained quadratic problem \eqref{eq:docp}, we apply the augmented Lagrangian method. To this end, let $\lambda \in \R^n$ and $\mu > 0$ and define the augmented Lagrangian function
\begin{align*}
    \mathcal{L}(z,\lambda;\mu) &:= \frac{1}{2}z^\top \mathcal Sz + \frac{\mu}{2} \|\mathcal Pz-x^0\|^2_2 -\lambda^\top (\mathcal Pz-x^0)\\
    &\phantom{:}= \frac12 z^\top(\mathcal{S} + \mu \mathcal{P}^\top \mathcal{P})z - (\mu x^0 + \lambda)^\top\mathcal{P}z \\&
    \qquad + \left(\frac{\mu}{2}x^0-\lambda\right)^\top x^0.
\end{align*}
Then, the unconstrained problem
\begin{align*}    \operatorname*{minimize}_z\  \mathcal{L}(z,\lambda_k;\mu_k)
\end{align*}
is successively solved for an increasing sequence of penalty parameters $\mu_k>0$ and corresponding suitable choice of the Lagrange multiplier $\lambda_k \in \R^n$. The algorithm is shown in Algorithm~\ref{alg:al}. 

\begin{algorithm}
\caption{Augmented Lagrangian method}\label{alg:al}
\begin{algorithmic}[1]
\Require $\mu_0>0$, $\mu_\Delta>0$, tolerance $\delta >0$, $z_0^s\in \R^{q}$, $\lambda_0\in \R^n$.
\Ensure approximate minimizer $z_*$ of \eqref{eq:docp}.
\State $k \gets 0$.
\While{not converged}
    \State $z_k \gets $ approx.\ min.\ of $\mathcal{L}(\cdot,\lambda_k;\mu_k)$ starting at $z_0^s$ satisfying $\|\nabla_z \mathcal{L}(z_k,\lambda_k;\mu_k)\| \leq \delta$.
    \State $\lambda_{k+1} \gets \lambda_k - \mu_k (\mathcal Pz-x^0)$.
    \State $\mu_{k+1} \gets \mu_k + \mu_\Delta$.
    \State $z_0^s \gets z_k$.
    \State $k \gets k+1$.
\EndWhile
\State $z_* \gets z_{k-1}$.
\end{algorithmic}
\end{algorithm}

Note that, in contrast to plain quadratic penalty methods (i.e., setting $\lambda_k = 0$) which often suffer from inherent ill-conditioning when increasing the penalty parameter$\mu_k$, see~\cite[Section 17.1]{NoceWrig06}, the augmented Lagrangian method includes a Lagrange multiplier $\lambda_k$ which alleviates this problem, cf.~\cite[Section 17.3]{NoceWrig06}.

\subsection{Inner loop: Limited memory BFGS}
We briefly derive the BFGS-method to solve the subproblem in Algorithm~\ref{alg:al} considered with the minimization of
\begin{align*}
     f_k(z) := \mathcal{L}(z,\lambda_k;\mu_k)
\end{align*}
over $z\in \R^q$.
Direct calculations yield that the gradient w.r.t.\ $z$ is given by
\begin{align*}
    \nabla f_k(z) &= \nabla_z \mathcal{L}(z,\lambda_k;\mu_k) \\ &=  (\mathcal{S} + \mu_k \mathcal{P}^\top \mathcal{P})z - \mathcal P^\top(\mu_k x^0 + \lambda_k).
\end{align*}
As the model is quadratic, and as standard in optimization, we can compute directly the minimizer for line searches, i.e., for the iterate $z_j$ and a search direction $p_j$
\begin{align*}
    \alpha_j = \operatorname*{argmin}_{\alpha \in \R} f_k(z_j + \alpha p_j) =  \frac{-\nabla f_k(z_j)^\top p_j}{p_j^\top (\mathcal S + \mu_k \mathcal P^\top \mathcal P)p_j}.
\end{align*}

The BFGS method is stated in Algorithm~\ref{alg:bfgs}. Whereas this algorithm is standard in numerical optimization and can be shown to converge superlinearly~\cite[Chapter 6]{NoceWrig06}, we state it here to highlight the fact that we may implement it without storing the matrices $\mathcal{S}$ and $\mathcal{P}$, i.e., by using the memory efficient DFT-based approach highlighted in Section~\ref{sec:dft}. 
\begin{algorithm}
\caption{BFGS method to solve $\min_z f_k(z)$}\label{alg:bfgs}
\begin{algorithmic}[1]
\Require $z_0\in \R^q$, $H_0 = I_{q}$, tolerance $\varepsilon>0$, parameters $\lambda_k,\,\mu_k$ defining $f_k = \mathcal L(\,\cdot\,, \lambda_k;\mu_k)$.
\Ensure approximate minimizer $z_*$ of $f_k$.
\State $j \gets 0$.
\While{$\|\nabla f_k(z_j)\| \geq \varepsilon$}
    \State $p_j \gets -H_j \nabla f_k(z_j)$. \Comment{Compute search direction.}
    \State $\alpha_j \gets \frac{-\nabla f_k(z_j)^\top p_j}{p_j^\top (\mathcal S + \mu_k \mathcal P^\top \mathcal P)p_j}$. \Comment{Exact line search.}
    \State $z_{j+1} \gets z_j + \alpha_j p_j$. \Comment{Update iterate.}
    \State $s_j \gets z_{j+1}-z_j$.
    \State $y_j \gets \nabla f_k(z_{j+1}) - \nabla f_k(z_j)$.
    \State $\rho_j \gets (s_j^\top y_j)^{-1}$.
    \State $H_{j+1} \gets \big(I-\rho_js_jy_j^\top\big)H_j\big(I-\rho_jy_js_j^\top\big) + \rho_js_js_j^\top $.
    \State $j\gets j+1$.
\EndWhile
\State $z_* \gets z_{j}$.
\end{algorithmic}
\end{algorithm}

In limited memory BFGS~\cite{LiuNoce89}, the update step for the inverse Hessian approximation $H_j$ is modified. Instead of including the full history of $\{s_0,\ldots,s_j\}$ and $\{y_0,\ldots,y_j\}$, only a window of fixed size of the most recent quantities is used.

\section{Numerical simulations}

We demonstrate the proposed algorithm for the quadratic problem~\eqref{eq:docp} with various state dimensions $n$, input dimensions $m$, and prediction horizons $L$. 
For each experiment, the LTI system $(A,B)$ is drawn randomly such that $A$ has spectral radius equal to $0.9$ to ensure stability of the system, and such that $A$ and $B$ are sparse with roughly half of their entries equal to zero. The input sequence $\tilde{\bm u}_N$ of the data trajectory $(\tilde{\bm x}_N, \tilde{\bm u}_N)$ is drawn randomly, where the signal length $N$ is optimized regarding the speedup of the DFT, cf.\ Remark~\ref{rm:fft}. Further, the control task~\eqref{eq:ocp} is extended to setpoint tracking, where the setpoint is some equilibrium of system \eqref{eq:LTI}. All calculations were performed in \texttt{MATLAB R2023b} on a \texttt{Apple M1 Pro} CPU. 

\begin{figure}[htb]
    \centering
    \input{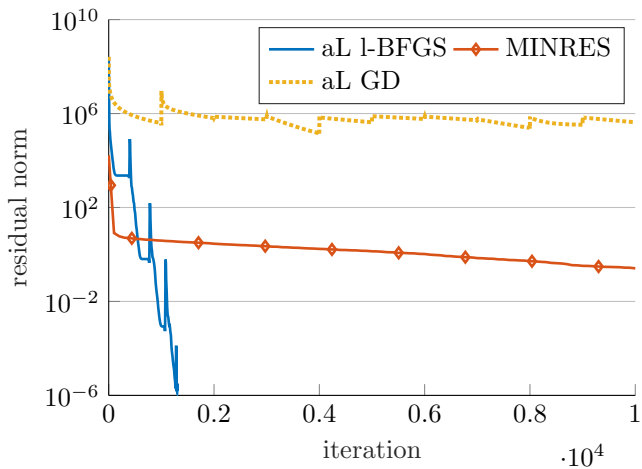}
    \caption{Residuals over iterations for the proposed method (aL lBFGS), MINRES, and gradient descent employing an augmented Lagrangian approach (aL GD) for $n=100$, $m=L = 50$.}
    \label{fig:residuum}
\end{figure}

In Figure~\ref{fig:residuum}, the progression of the residuals regarding the iteration number is shown for the proposed method (aL lBFGS), MINRES~\cite{paige1975solution}, and gradient descent combined with an augmented Lagrangian (aL GD) for the parameters $n=100$, $m=L = 50$. The residual $r_k$, which is defined as
\begin{equation*}
    r_k := \begin{bmatrix}
        \mathcal S z_k - \mathcal P^\top \lambda_k\\
        \mathcal P z_k - x^0
    \end{bmatrix}
\end{equation*}
with approximated primal and dual solutions $z_k$ and $\lambda_k$ at iteration $k$, indicates the violation of the first-order Karush--Kuhn--Tucker condition for the constrained problem~\eqref{eq:docp}. It is evident that neither gradient descent nor MINRES perform well. 

While running, the BFGS algorithm constructs a preconditioner $\mathcal B$ of the matrix $\mathcal S$ which is improved in every iteration step. In Figure~\ref{fig:condition} the condition numbers\footnote{The condition number is defined as the ratio between the largest and smallest non-zero singular value.} $\kappa(\mathcal{S})$ and $\kappa(\mathcal B\mathcal{S})$ of $\mathcal S$ and $\mathcal B\mathcal S$ are depicted in the case $m=L=50$. Even though the  condition drastically improves during a BFGS run, $\kappa(\mathcal B\mathcal S)$ seems to grow similarly to $\kappa(\mathcal S)$ in the state dimension.

\begin{figure}[htb]
    \centering
\definecolor{mycolor1}{rgb}{0.00000,0.44700,0.74100}%
\definecolor{mycolor2}{rgb}{0.85000,0.32500,0.09800}%
\begin{tikzpicture}

\begin{axis}[%
width=7cm,
height=5cm,
at={(1.011in,0.642in)},
scale only axis,
xmode=log,
xmin=10,
xmax=500,
xminorticks=true,
xlabel style={font=\color{white!15!black}},
xlabel={state dimension $n$},
ymode=log,
ymin=100,
ymax=100000000000,
yminorticks=true,
ylabel style={font=\color{white!15!black}},
ylabel={condition number},
axis background/.style={fill=white},
axis x line*=bottom,
axis y line*=left,
%xmajorgrids,
%xminorgrids,
ymajorgrids,
%yminorgrids,
legend style={at={(0.03,0.97)}, anchor=north west, legend cell align=left, align=left, draw=white!15!black}
]
\addplot [color=mycolor1, line width=1.0pt, mark=o, mark options={solid, mycolor1}]
  table[row sep=crcr]{%
10	23498596.4034423\\
20	4525157.19239954\\
30	1935531.58734895\\
40	641296.176514722\\
50	1602781.58172546\\
60	1745843.56211101\\
70	1449499.52694169\\
80	3900384.13322225\\
100	494256.505972761\\
120	2143788.33951515\\
140	2926454.66683804\\
160	5176128.59073883\\
180	23164014.4558788\\
200	22769909.780537\\
220	31747323.9008804\\
240	24013335.8508286\\
260	90157086.7608142\\
280	274359123.596429\\
300	231620996.953582\\
350	967219908.770036\\
400	2469765569.89702\\
450	14702217521.3459\\
500	47450536170.4624\\
};
\addlegendentry{$\kappa(\mathcal S)$}

\addplot [color=mycolor2, line width=1.0pt, mark=diamond, mark options={solid, mycolor2}]
  table[row sep=crcr]{%
10	729.034875960414\\
20	689.638739452968\\
30	616.41607681864\\
40	488.22692428402\\
50	492.528056819398\\
60	480.612637642077\\
70	390.133348128199\\
80	404.023387738628\\
100	382.628263238687\\
120	397.252741415751\\
140	589.754916799367\\
160	901.357221825878\\
180	1961.2347453505\\
200	1270.6978508681\\
220	3337.97782062306\\
240	4296.88737323942\\
260	3040.25360446489\\
280	11577.1314250715\\
300	6584.35283808449\\
350	26781.3315521859\\
400	161202.069629635\\
450	696723.110493103\\
500	70038427.5362785\\
};
\addlegendentry{$\kappa(\mathcal B\mathcal S)$}

\end{axis}
\end{tikzpicture}%
    \caption{The condition number of the matrix $\mathcal S$ (blue circles) and $\mathcal B\mathcal S$ (red diamonds) for increasing state dimension $n$ and fixed values $m = L = 50$.
    }
    \label{fig:condition}
\end{figure}
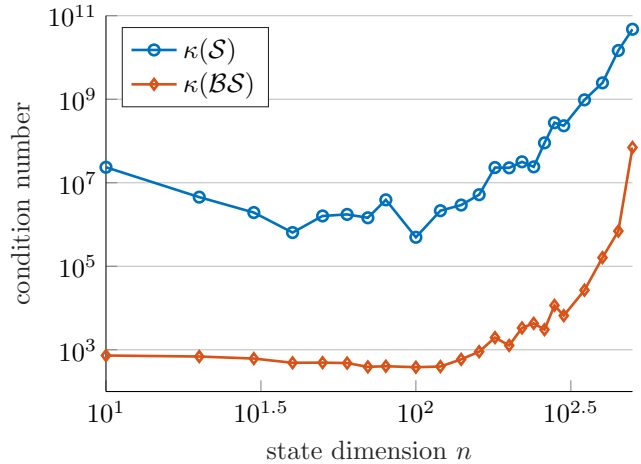

\begin{figure}[htb]
    \centering
\begin{tikzpicture}[scale=1.0]

\begin{axis}[%
width=7cm,
height=5cm,
at={(1.011in,0.642in)},
scale only axis,
xmode=log,
xmin=10,
xmax=1000,
xtick={  10,  100, 1000},
%xminorticks=true,
xlabel style={font=\color{white!15!black}},
xlabel={state dimension $n$},
ymode=log,
ymin=1000,
ymax=50000,
%yminorticks=true,
ylabel style={font=\color{white!15!black}},
ylabel={total number of BFGS iterations},
axis background/.style={fill=white},
axis x line*=bottom,
axis y line*=left,
%xmajorgrids,
%xminorgrids,
ymajorgrids,
%yminorgrids,
legend style={at={(0.03,0.97)}, anchor=north west, legend cell align=left, align=left, draw=white!15!black}
]
\addplot [color=mycolor1, line width=1.0pt, mark=o, mark options={solid, mycolor1}]
  table[row sep=crcr]{%
10	1390\\
20	1326\\
30	1370\\
40	1561\\
50	1402\\
60	1447\\
70	1648\\
80	1780\\
90	1769\\
100	2059\\
120	1925\\
140	2048\\
160	1938\\
180	2057\\
200	2331\\
220	2272\\
240	2521\\
260	2689\\
280	2849\\
300	2983\\
350	4097\\
400	5015\\
450	6132\\
500	7736\\
600	11440\\
650	14837\\
700	16024\\
750	19172\\
800	22440\\
850	28143\\
900	32074\\
};
\addlegendentry{$m=L=50$}

\addplot [color=mycolor2, line width=1.0pt, mark=diamond, mark options={solid, mycolor2}]
  table[row sep=crcr]{%
10	2777\\
20	2699\\
30	2722\\
40	2917\\
50	2450\\
60	2935\\
70	3057\\
80	2903\\
90	3050\\
100	3053\\
120	3145\\
140	3106\\
160	3020\\
180	3062\\
200	3251\\
220	3477\\
240	3450\\
260	3250\\
280	3135\\
300	4161\\
350	4408\\
400	4467\\
450	4431\\
500	4848\\
550	5087\\
600	5596\\
650	6386\\
700	7481\\
750	8584\\
800	9078\\
850	10889\\
900	11703\\
950	13152\\
1000	15252\\
};
\addlegendentry{$m=L=100$}

% \addplot [color=mycolor3]
%   table[row sep=crcr]{%
% 10	0.2\\
% 20	2.08164799306237\\
% 30	7.97645477548618\\
% 40	20.5063678889979\\
% 50	42.4742501084005\\
% 60	76.8161340165734\\
% 70	126.573725544978\\
% 80	194.876414667975\\
% 90	284.928557876254\\
% 100	400\\
% 120	718.565038634059\\
% 140	1177.79506598022\\
% 160	1805.61508979173\\
% 180	2630.5498499525\\
% 200	3681.64799306237\\
% 220	4988.42334107897\\
% 240	6580.80804108425\\
% 260	8489.11431278702\\
% 280	10744.0026208049\\
% 300	13376.4547754862\\
% 350	21815.3834803036\\
% 400	33306.3678889979\\
% 450	48354.7980635556\\
% 500	67474.2501084005\\
% 600	120016.134016573\\
% 650	154499.266113609\\
% 700	195173.725544978\\
% 750	242583.294098675\\
% 800	297276.414667975\\
% 850	359805.879550858\\
% 900	430728.557876254\\
% };
% \addlegendentry{$\mathcal O(n^3\log n)$}

\end{axis}
\end{tikzpicture}%
    \caption{Total number of BFGS iterations in dependence of the state dimension $n$ in \eqref{eq:docp} for fixed $m=L=50$ (blue circles) and for $m=L=100$ (red diamonds).}
    \label{fig:iterations}
\end{figure}
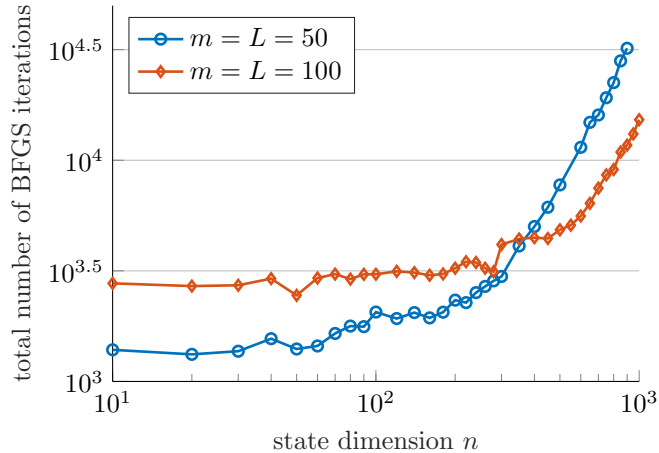

\begin{figure}[htb]
    \centering
\begin{tikzpicture}

\begin{axis}[%
width=7cm,
height=5cm,
at={(1.011in,0.642in)},
scale only axis,
xmode=log,
xmin=10,
xmax=1000,
xminorticks=true,
xlabel style={font=\color{white!15!black}},
xlabel={state dimension $n$},
ymode=log,
ymin=1,
ymax=40000,
yminorticks=true,
ylabel style={font=\color{white!15!black}},
ylabel={total execution time [s]},
axis background/.style={fill=white},
axis x line*=bottom,
axis y line*=left,
%xmajorgrids,
%xminorgrids,
ymajorgrids,
%yminorgrids,
legend style={at={(0.03,0.97)}, anchor=north west, legend cell align=left, align=left, draw=white!15!black}
]
\addplot [color=mycolor1, line width=1.0pt, mark=o, mark options={solid, mycolor1}]
  table[row sep=crcr]{%
10	7.58738754166667\\
20	8.23186425\\
30	10.3835744166667\\
40	15.803551125\\
50	16.7540311666667\\
60	18.9631968333333\\
70	25.029640625\\
80	29.4905742916667\\
90	32.445343875\\
100	36.5514541666667\\
120	45.3822892083333\\
140	58.7641657916667\\
160	95.7679045416667\\
180	118.575747916667\\
200	162.7705475\\
220	173.507156416667\\
240	213.386521208333\\
260	261.114360791667\\
280	322.444013833333\\
300	363.181651916667\\
350	641.65055825\\
400	1152.40540629167\\
450	1837.70954908333\\
500	2722.04920670833\\
600	5855.59998395833\\
650	9065.40101033333\\
700	10736.349513875\\
750	13742.9423680417\\
800	19419.9082079167\\
850	36886.49276975\\
900	59074.7059302083\\
};
\addlegendentry{$m=L=50$}

\addplot [color=mycolor2, line width=1.0pt, mark=diamond, mark options={solid, mycolor2}]
  table[row sep=crcr]{%
10	108.346701458333\\
20	110.616583583333\\
30	116.958738166667\\
40	147.774229958333\\
50	132.896564\\
60	156.692250416667\\
70	194.605218916667\\
80	210.342224375\\
90	248.796162666667\\
100	284.125464958333\\
120	386.27991775\\
140	421.910869125\\
160	416.153126166667\\
180	485.998371916667\\
200	702.252955625\\
220	1297.84830279167\\
240	1158.623336\\
260	1148.83503870833\\
280	1090.953252\\
300	1656.62454666667\\
350	1884.37739879167\\
400	1843.42203891667\\
450	2072.88753325\\
500	2627.21164408333\\
550	3247.76546920833\\
600	4067.95125758333\\
650	5438.44241483333\\
700	6913.97023395833\\
750	8725.25935641667\\
800	11007.932636375\\
850	14346.26860025\\
900	18376.167983125\\
950	23692.9684314167\\
1000	38681.712511875\\
};
\addlegendentry{$m=L=100$}

% \addplot [color=mycolor3]
%   table[row sep=crcr]{%
% 10	1e-06\\
% 20	8.32659197224948e-05\\
% 30	0.00107682139469063\\
% 40	0.00656203772447933\\
% 50	0.0265464063177503\\
% 60	0.0829614247378993\\
% 70	0.217073939309637\\
% 80	0.498883621550016\\
% 90	1.03856459345894\\
% 100	2\\
% 120	6.20840193379827\\
% 140	16.1593483052486\\
% 160	36.9789970389347\\
% 180	76.7068336246148\\
% 200	147.265919722495\\
% 220	265.583658679044\\
% 240	454.865451799743\\
% 260	746.023365807723\\
% 280	1179.26172765954\\
% 300	1805.82139469063\\
% 350	4676.67283359009\\
% 400	10658.0377244793\\
% 450	22031.6548677075\\
% 500	42171.4063177503\\
% 600	129617.424737899\\
% 650	212146.804782249\\
% 700	334722.939309637\\
% 750	511699.135989392\\
% 800	761027.621550016\\
% 850	1104828.92889585\\
% 900	1570005.59345894\\
% };
%\addlegendentry{$\mathcal O(n^6\log n)$}

\end{axis}
\end{tikzpicture}%
    \caption{Total execution time in dependence of the state dimension $n$ in \eqref{eq:docp} for fixed $m=L=50$ (blue circles) and for $m=L=100$ (red diamonds).}
    \label{fig:exec_time}
\end{figure}
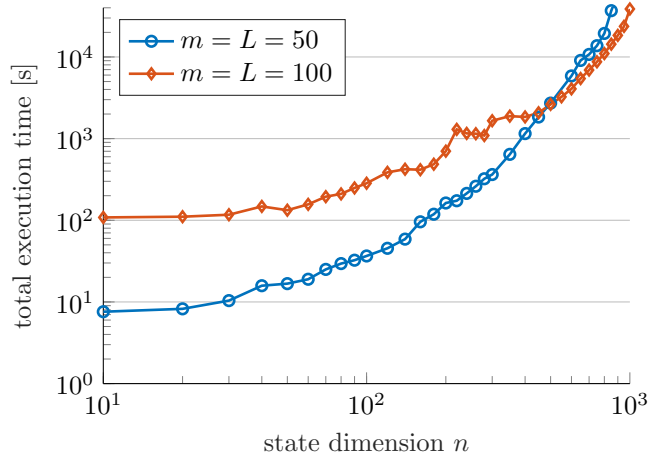

\begin{figure}[htb]
    \centering
\begin{tikzpicture}

\begin{axis}[%
width=7cm,
height=5cm,
at={(1.011in,0.642in)},
scale only axis,
xmode=log,
xmin=10,
xmax=1000,
xminorticks=true,
xlabel style={font=\color{white!15!black}},
xlabel={state dimension $n$},
ymode=log,
ymin=0.001,
ymax=3,
yminorticks=true,
ylabel style={font=\color{white!15!black}},
ylabel={mean execution time per BFGS iteration [s]},
axis background/.style={fill=white},
axis x line*=bottom,
axis y line*=left,
%xmajorgrids,
%xminorgrids,
ymajorgrids,
%yminorgrids,
legend style={at={(0.03,0.97)}, anchor=north west, legend cell align=left, align=left, draw=white!15!black}
]
\addplot [color=mycolor1, line width=1.0pt, mark=o, mark options={solid, mycolor1}]
  table[row sep=crcr]{%
10	0.0054585521882494\\
20	0.00620804242081448\\
30	0.00757925139902676\\
40	0.010123991752082\\
50	0.0119500935568236\\
60	0.0131051809490901\\
70	0.0151878887287621\\
80	0.0165677383661049\\
90	0.018341064937818\\
100	0.0177520418487939\\
120	0.0235752151731602\\
140	0.0286934403279622\\
160	0.0494158434167527\\
180	0.0576449916950251\\
200	0.0698286347061347\\
220	0.0763675864509977\\
240	0.0846436022246463\\
260	0.0971046339872319\\
280	0.113177962033462\\
300	0.121750469968712\\
350	0.156614732304125\\
400	0.229791706139914\\
450	0.29969170728691\\
500	0.351867787837168\\
600	0.511853145450903\\
650	0.610999596302038\\
700	0.670016819388105\\
750	0.716823616108996\\
800	0.865414804274361\\
850	1.31068090714387\\
900	1.84182533922206\\
};
\addlegendentry{$m=L=50$}

\addplot [color=mycolor2, line width=1.0pt, mark=diamond, mark options={solid, mycolor2}]
  table[row sep=crcr]{%
10	0.0390157369313408\\
20	0.0409842843954551\\
30	0.0429679420156747\\
40	0.0506596605959319\\
50	0.0542434955102041\\
60	0.0533874788472459\\
70	0.0636588874441173\\
80	0.0724568461505339\\
90	0.0815725123497268\\
100	0.0930643514439349\\
120	0.122823503259141\\
140	0.135837369325499\\
160	0.137799048399559\\
180	0.158719259280427\\
200	0.216011367463857\\
220	0.373266696229988\\
240	0.335832851014493\\
260	0.353487704217949\\
280	0.347991467942584\\
300	0.398131349835777\\
350	0.427490335479053\\
400	0.412675629934333\\
450	0.467814834856692\\
500	0.541916593251513\\
550	0.638444165364327\\
600	0.726939109646771\\
650	0.851619545072555\\
700	0.924204014698347\\
750	1.01645612260213\\
800	1.21259447415455\\
850	1.31750101940031\\
900	1.57021003017389\\
950	1.80147266053959\\
1000	2.53617312561467\\
};
\addlegendentry{$m=L=100$}

\addplot [color=mycolor3]
  table[row sep=crcr]{%
10	0.000204725443947948\\
20	0.00106541577380762\\
30	0.0027216387417369\\
40	0.00524771908729382\\
50	0.00869555970979848\\
60	0.0131051809490901\\
70	0.018509187253101\\
80	0.0249351003174285\\
90	0.0324067363914804\\
100	0.0409450887895897\\
120	0.0612952277249302\\
140	0.0861159349151797\\
160	0.115517297142727\\
180	0.149594579405204\\
200	0.188431754959942\\
220	0.232103905250614\\
240	0.280678926614\\
260	0.334218787771193\\
280	0.392780483271821\\
300	0.456416773726844\\
350	0.638023438271414\\
400	0.852332619046099\\
450	1.09993972228609\\
500	1.38136958084972\\
600	2.04752969312162\\
650	2.4330716034416\\
700	2.85407340065505\\
750	3.31086483858046\\
800	3.80375287300972\\
850	4.33302461330185\\
900	4.89894973512642\\
};
\addlegendentry{$\mathcal O(n^2\log n)$}

\end{axis}
\end{tikzpicture}%
   \caption{Mean execution time per BFGS iteration (total execution time divided by the total 
  number of iterations divided) for increasing $n$ and fixed $m,L$ with $m=L=50$ (blue circles) and $m=L=100$ (red diamonds).}
    \label{fig:mean_exec_time}
\end{figure}
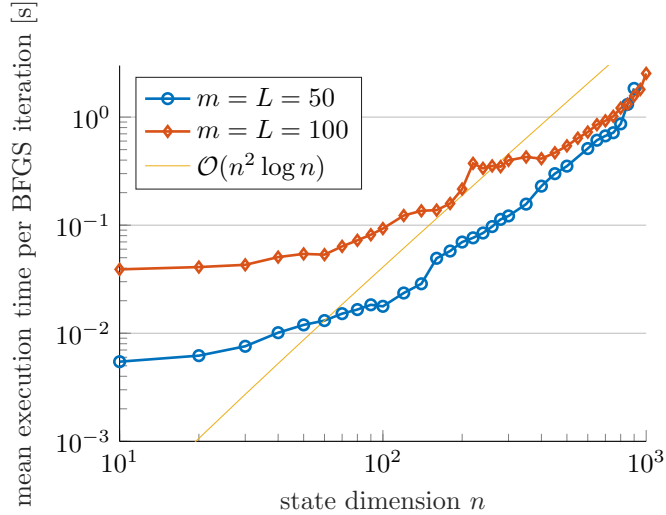

In Figures~\ref{fig:iterations} and~\ref{fig:exec_time}, the total number of BFGS iterations and the total execution time for solving \eqref{eq:docp} with the proposed method are shown. The connection between Figures~\ref{fig:iterations} and~\ref{fig:exec_time} can be seen in Figure~\ref{fig:mean_exec_time}, which illustrates the mean execution time per BFGS iteration given as the ratio of total execution time and total number of BFGS iterations. The theoretical asymptotical behavior of $\cO(n^2\log n)$ as $n\to\infty$ for the DFT-based matrix-vector multiplication is reflected in Figure~\ref{fig:mean_exec_time}, cf.\ Table~\ref{tab:2}.

\section{Conclusion}
We proposed an iterative algorithm to solve data-driven optimal control problems involving Hankel-like matrices. 
To this end, we tackled both the high dimension and bad conditioning of the problem by means of a DFT-based factorization of the respective  
matrices (enabling fast matrix-vector multiplication) in combination 
with an lBFGS method. 
The latter allows to construct a preconditioner iteratively, without relying on  storing the matrix, which allows for the solution of problems with high dimension stemming from many data samples.

Future work may exploit recently-proposed concepts by~\cite{berberich2023quantitative,coulson2022quantitative} in combination with the proposed algorithmic advancements or the application of the derived results w.r.t.\ energy systems based on our previous work~\cite{schmitz2022data}.

\section*{Code Availability}
The code that produced the numerical results is available at
\begin{center}
    \url{https://github.com/schmitzph/FastDeePC}.
\end{center}

\section*{Aknowledgement}
Philipp Schmitz is grateful for the support from the Carl Zeiss Foundation (VerneDCt, Project No.\ 2011640173). Karl~Worthmann gratefully acknowledges funding by the German Research Foundation (DFG, Project-ID 507037103).

\bibliographystyle{abbrv}
\bibliography{references}

\end{document}